\documentclass[11pt]{article}
\usepackage{amssymb,latexsym, amsmath}
\usepackage{graphicx}

\usepackage{times}
\usepackage[T1]{fontenc}
\usepackage{mathrsfs}
\usepackage{epsfig}
\usepackage{color}

\newtheorem{theorem}{Theorem}[section]

\newtheorem{defi}[theorem]{Definition}
\newtheorem{lemma}[theorem]{Lemma}

\newtheorem{cj}[theorem]{Conjecture}

\def\slfrac#1#2{\hbox{\kern.1em %
 \raise.5ex\hbox{\the\scriptfont0 #1}\kern-.11em %
 /\kern-.15em\lower.25ex\hbox{\the\scriptfont0 #2}}}

\newcommand{\eqn}[1]{(\ref{#1})}

\newcommand{\eeq}{\end{equation}}
\newcommand{\beql}[1]{\begin{equation}\label{#1}}

\newcommand{\loc}{loc}

\newcommand{\RR}{{\mathbb R}}

\newcommand{\NN}{{\mathbb N}}

\newcommand{\dM}{{N}^{1}}
\newcommand{\dN}{{D}}

\newcommand{\sB}{{\cal B}}
\newcommand{\sBp}{{\cal B}^{'}}

\newcommand{\sBo}{{B_{\emptyset}}}

\newcommand{\sG}{{\cal G}}
\newcommand{\sH}{{\cal H}}

\newcommand{\sP}{{\cal P}}

\newcommand{\tauL}{\tau^{L}}
\newcommand{\tauS}{\tau^{S}}

\newcommand{\ddL}{\Omega^{L}}

\newcommand{\muT}{{\mu_S}}
\DeclareMathOperator{\meas}{{meas}}

\makeatletter
\def\@sect#1#2#3#4#5#6[#7]#8{\ifnum #2>\c@secnumdepth
     \def\@svsec{}\else
     \refstepcounter{#1}\edef\@svsec{\csname the#1\endcsname.\hskip .75em }\fi
     \@tempskipa #5\relax
      \ifdim \@tempskipa>\z@
        \begingroup #6\relax
          \@hangfrom{\hskip #3\relax\@svsec}{\interlinepenalty \@M #8\par}%
        \endgroup
       \csname #1mark\endcsname{#7}\addcontentsline
         {toc}{#1}{\ifnum #2>\c@secnumdepth \else
                      \protect\numberline{\csname the#1\endcsname}\fi
                    #7}\else
        \def\@svsechd{#6\hskip #3\@svsec #8\csname #1mark\endcsname
                      {#7}\addcontentsline
                           {toc}{#1}{\ifnum #2>\c@secnumdepth \else
                             \protect\numberline{\csname the#1\endcsname}\fi
                       #7}}\fi
     \@xsect{#5}}
\def\@begintheorem#1#2{\it \trivlist \item[\hskip \labelsep{\bf #1\ #2.}]}

\def\plain{plain}\ifx\fmtname\plain\csname fi\endcsname
     
     \input docstrip
     \preamble

     Do not distribute the stripped version of this file.
     The checksum in the header refers to the documented version.

     \endpreamble
     \generateFile{here.sty}{t}{\from{here.doc}{}}
     \endinput
\fi
\ifcat a\noexpand @\let\next\relax\else\def\next{%
    \documentstyle[here,doc]{article}\MakePercentIgnore}\fi\next
\ifx\@Hxfloat\@Hundef\else\expandafter\endinput\fi
\let\@Hxfloat\@xfloat
\def\@xfloat#1[{\@ifnextchar{H}{\@HHfloat{#1}[}{\@Hxfloat{#1}[}}
\def\@HHfloat#1[H]{%
\expandafter\let\csname end#1\endcsname\end@Hfloat
\vskip\intextsep\vbox\bgroup\def\@captype{#1}\parindent\z@
\ignorespaces}
\def\end@Hfloat{\egroup\vskip \intextsep}

\makeatother

\begin{document}

\begin{center}
{\Large 
{\bf   Level Sets of the Takagi Function: Generic Level Sets} 
}\\

\vspace{1.5\baselineskip}
{\em Jeffrey C. Lagarias
\footnote{This author's work was supported by NSF Grants DMS-0500555 and DMS-0801029.} }\\
\vspace*{.2\baselineskip}
Dept. of Mathematics \\
University of Michigan \\
Ann Arbor, MI 48109-1043\\
\vspace*{1.5\baselineskip}
{\em Zachary Maddock}
\footnote{This author's work was supported by the NSF
  through a Graduate Research Fellowship.} \\
\vspace*{.2\baselineskip}
Dept. of Mathematics \\
Columbia University  \\
New York, NY  10027\\

\vspace*{2\baselineskip}
(July 26, 2011) \\
\vspace{3\baselineskip}
{\bf ABSTRACT}
\end{center}
The Takagi function $\tau: [0,1] \to [0, 1]$ 
is a 
continuous non-differentiable function
constructed by Takagi in 1903.
This paper studies the level sets $L(y)=\{ x: \tau(x)=y \}$ 
 of the Takagi function $\tau(x)$. It shows that for a  full 
 Lebesgue measure set of ordinates $y$,
 these level sets are  finite sets, but whose expected number of
 points  is infinite. 
 Complementing this, it  shows
 that the set of ordinates $y$ whose level set has 
  positive Hausdorff dimension is itself a set of
  full Hausdorff dimension $1$ (but Lebesgue measure zero).
  Finally it shows that the level sets have a nontrivial
Hausdorff dimension spectrum.
 The results are obtained using a notion of ``local level set"  introduced in 
a previous paper, 
 along with a singular measure parameterizing  such sets.

%
%
%
%
\setlength{\baselineskip}{1.0\baselineskip}

\section{Introduction}

The Takagi function $\tau(x)$ is  a function
defined on the unit interval $ x\in [0,1]$
which was introduced by Takagi \cite{Takagi} in 1903
as an example of a continuous nondifferentiable function.
It can be defined by 
\begin{equation}\label{eq101}
\tau(x) := \sum_{n=0}^{\infty} \frac{\ll 2^n x \gg}{2^n}
\end{equation}
where $\ll x \gg$ is the distance from $x$ to the nearest integer.
It has appeared  in a wide variety of contexts in probability theory, 
number theory and analysis,
including
Bernoulli convolutions \cite[p. 195]{HY84}, the distribution
of binary digit sums
( (\cite{Mir49}, \cite{Trollope68}, \cite{Delange75}),
and as fractals and functions satisfying self-similar
analogues of the Laplace equation (\cite{YHK97}). \smallskip
%
%

\begin{figure}[h]
  \begin{center}$                                                               
    \begin{array}{cc}
      \includegraphics[height=2in]{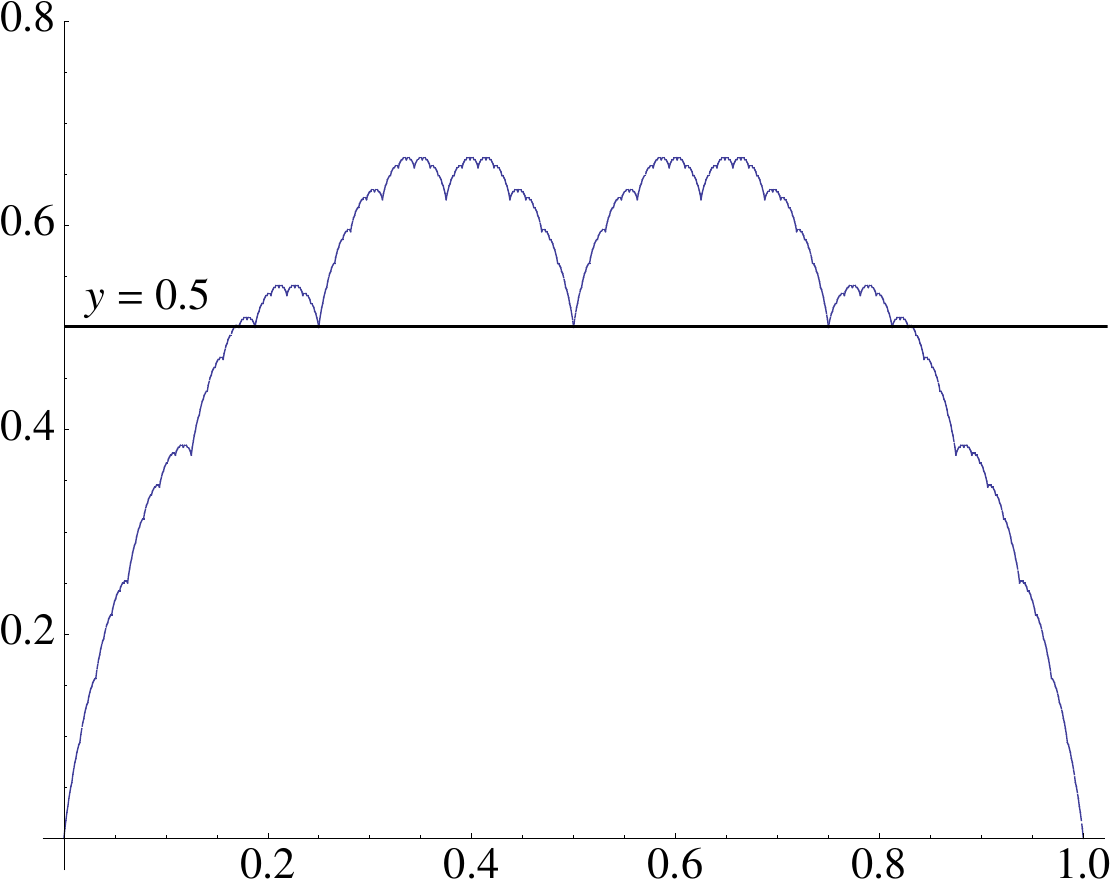} &
      \includegraphics[height=2in]{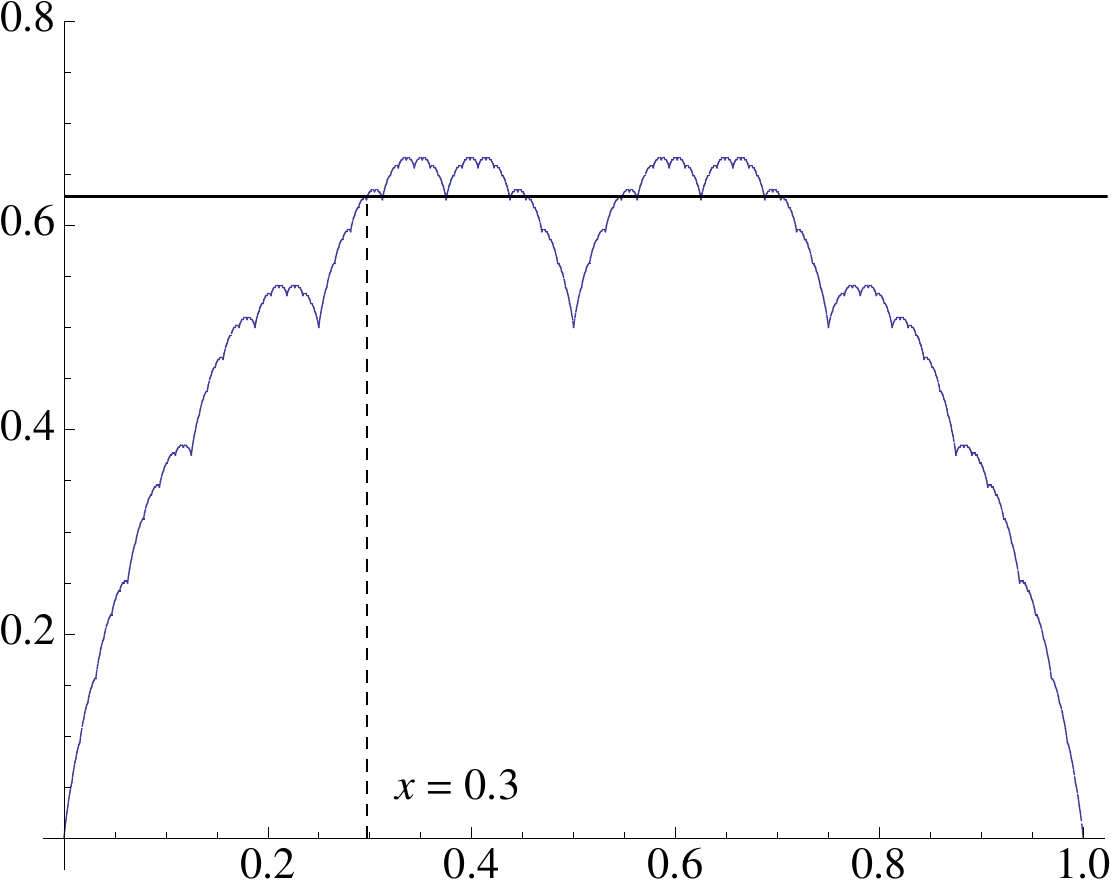}
    \end{array}$
  \end{center}
\caption{Ordinate level  set $L(y)$ at  $y=0.5$ and  abscissa level set
  $L(\tau(x))$ at   $x=0.3$.}
\label{fig11}
\end{figure}

This paper considers certain properties of the graph of the Takagi function
$$
\sG(\tau):= \{ (x, \tau(x)): ~0 \le x \le 1 \},
$$
which is pictured in Figure~\ref{fig11}. 
It is well known that the values of the Takagi function satisfy $0 \le \tau(x) \le \frac{2}{3}.$
It is also known that this graph has Hausdorff dimension $1$ in $\RR^2$, see Mauldin
and Williams \cite[Theorem 7]{MW86}). 
They add the remark  that they do not know whether the
$1$-dimensional Hausdorff measure of this graph is  finite or infinite  (\cite[p. 800]{MW86}). 
Here we  study  the structure of the level sets of this graph. We make the
following definition, which contains a special convention concerning
dyadic rationals which simplifies theorem statements. 

%
%

\begin{defi}~\label{de11}
{\em 
For $0 \le y \le \frac{2}{3}$  the 
{\em (global) level set}  $L(y)$ at level $y$ is 
$$
L(y):= \{ x: ~ \tau(x) =y, ~~0 \le x \le 1 \}. 
$$
We make the convention that 
$x$ specifies a particular binary expansion; so each dyadic rational
value  $x= \frac{m}{2^n}$  in a level set will appear twice, labeled by each of its two possible
 binary expansions. 
}
\end{defi}

Level sets have a complicated and interesting structure, depending on the value of $y$.
It is known that there are different levels $y$ where the level set  $L(y)$ 
is finite, countably infinite, or
uncountably infinite, respectively.
Concerning the size of level sets, measured by Hausdorff dimension, in 1984
Baba \cite{Baba84} showed that the level set $L(\frac{2}{3})$ has Hausdorff dimension
$\frac{1}{2}$, so is uncountable. 
The second author showed (\cite{Mad10}) that the Hausdorff
dimension of any level set is at most $0.699$ and conjectured 
Baba's example achieves the largest possible dimension.  
Finally,  relevant to our results here, in 2008 Buczolich \cite{Buz08} 
proved that, in the sense of Lebesgue measure, almost all level sets $L(y)$ are finite sets.\smallskip

The object of this paper is to   study properties of 
  ``generic" level sets of the Takagi function.
The  {\em ordinate (y-axis) notion of
genericity} (in the Lebesgue measure sense)  is  to draw an ordinate
$y$ at random using Lebesgue measure in $[0,  \frac{2}{3}]$, and 
 ask for properties of such level sets that hold for a full
Lebesgue measure of such sets $L(y)$.
The  {\em abscissa (x-axis) notion of genericity} is to draw a number $x$ at
random in $[0,1]$ with respect to Lebesgue measure, and then to ask 
for properties of level sets $L(\tau(x))$
that hold for a full  Lebesgue measure set of $x$. 
These two sampling methods of drawing level sets are pictured in
Figure \ref{fig11} above. We focus on ordinate notions of genericity,
but  also treat abscissa notions in order to obtain our results: 
we observe that these two notions of genericity give quite
different answers concerning  the structure of level sets.\smallskip

A   weaker notion
of ``generic set" is  ask for properties of
level sets that are  {\em generic in  the Hausdorff dimension
sense}, which means to  hold for  a set of (full) Hausdorff dimension $1$,
in the appropriate variable. 
We consider this notion in both  the  abscissa sense 
and the  ordinate sense. 
  \smallskip


\subsection{Generic level sets: Results} \label{sec11}

Our results are  obtained using methods  from our paper
\cite{LM10a}, which 
 introduced the notion of ``local level set" of the Takagi function.
We summarize relevant results in Sect. \ref{sec2} and \ref{sec3}.
 Local level sets are sets $L_x^{loc}$ determined
locally by combinatorial operations on the binary expansion of a real number $x$; 
they are closed sets and each level set  is partitioned into a disjoint union of local level sets. 
The structure
of local level sets $L_x^{loc}$ is completely analyzable: they are either finite sets or Cantor
sets. Information of the  Hausdorff dimension of such sets can
 be deduced from properties of the binary expansion of $x$. 
 In \cite{LM10a} we introduced a certain subset $\ddL \subset [0,1]$,
called the {\em deficient digit set},  made up of 
the left endpoints of all local level sets; this set parameterizes all local level
sets.  We  showed 
 that $\ddL$  is a closed set of Lebesgue measure zero. 
We also introduced a new singular function $\mu_S$, 
the {\em Takagi singular function}, whose points of increase are  supported on $\ddL$.
This function was  used to show that 
on almost all  levels $y$  there 
are finitely many local level sets; 
however there is also a dense set of levels in $[0, \frac{2}{3}]$
which have infinitely many local level sets. 
In  Theorem \ref{th34} below  we introduce
 singular measure, the Takagi singular measure, associated to
 the Takagi singular function,
 which plays an important role in this paper.
\medskip

 Our first  result concerns  
 the cardinality  of  ``generic"   ordinate level sets in the Lebesgue measure sense.
 
%
%
\begin{theorem}~\label{th12}{\em (Ordinate generic level sets)}
 (1) For a full Lebesgue measure set of ordinate points $y \in  [0, \frac{2}{3}]$ the
 level set $L(y)$ is a finite set. 
 
  (2) For  a random level set $L(y)$ with level $y$ drawn 
  uniformly from $y \in [0, \frac{2}{3}]$, the expected number of elements in $L(y)$  is infinite.
\end{theorem} 
 
 We prove Theorem~\ref{th12}  in Sect. \ref{sec4} and Sect. \ref{sec5}.
 This result   is proved using explicit calculations of the Takagi
 singular measure of various subsets of $\ddL$ given in the fine decomposition
 of the deficient digit set $\ddL$ made in Sect. \ref{sec4} below.
 These calculations make use
 of  self-similarity properties of the Takagi singular measure.
\medskip

Part  (1) of this result was first
 proved  in 2008 by Buczolich \cite{Buz08}.
He  proves the almost everywhere finiteness of level sets
by a  method  that directly studies the graph of the Takagi function.
His proof shows the graph $\sG(\tau)= \sG_{I} \bigcup \sG_{R}$ 
(nonconstructively) partitions into an irregular
$1$-set $\sG_{R}$ and a regular $1$-set $\sG_{R}$, and that the irregular set $\sG_{I}$ has
$y$-projection of Lebesgue measure $0$ and $x$-projection of full measure $1$.
Here an {\em irregular $1$-set} or {\em purely unrectifiable $1$-set} is a set in $\RR^2$
of Hausdorff dimension $1$ that intersects every continuously differentiable
curve in a set of $\sH^{1}$-measure zero.
By Besicovich's theorem such  a set has $1$-dimensional projections
of measure $0$ in almost all directions,  see Falconer \cite[Theorem 6.1.3]{Fa85}.
A {\em regular $1$-set} is a set that can be covered by countably many 
rectifiable curves. 
In comparison our
proof of Theorem~\ref{th12} uses the Takagi singular measure, whose support lies in $\ddL$.
The part of the graph $\sG(\tau)$ that lies above $\ddL$ is covered by a single
rectifiable curve, the flattened Takagi function described in \cite[Sect. 5]{LM10a}
(which is of bounded variation  \cite[Theorem  5.4]{LM10a}), so it  should belong to the 
regular part $\sG_{R}$ of Buczolich's partition. \smallskip

 We contrast Theorem \ref{th12} with what is known about ``generic" abscissa  level
 sets in the Lebesgue measure sense.
In \cite[Theorem 1.4]{LM10a} we showed that a ``generic" 
local level set $L_x^{loc}$ 
obtained by drawing $x$ with the uniform distribution on $[0,1]$  (Lebesgue
measure) is with probability one an uncountable set of Hausdorff dimension $0$.
An immediate consequence is: {\em  For a full Lebesgue
measure set of abscissa points $x \in [0, 1]$ the level set $L(\tau(x))$ is uncountable.}
Here we advance  the following conjecture for abscissa level sets,
concerning their Hausdorff dimension.

%
%
\begin{cj}~\label{cg12}
{\rm (Abscissa generic level sets)}
  A full Lebesgue measure set of abscissa  points $x \in  [0, 1]$ have
 level sets $L(\tau(x))$ that  are uncountable and have Hausdorff dimension $0$.
\end{cj}

We also prove  generic results in the weaker Hausdorff dimension sense.
Let $\dim_{H}(\Gamma)$ denote the Hausdorff dimension of  a set $\Gamma$.
Our next two results show that the set of level sets that are large in the sense
of having positive Hausdorff dimension are themselves  generic in the 
Hausdorff dimension sense. We prove results in both the abscissa generic
case and the ordinate generic case. 
%
%
%
\begin{theorem}~\label{th13a} {\em (Positive Hausdorff dimension abscissa level sets)}
Let $\Gamma_H^{abs}$ be  set of abscissas $x \in [0,1]$ such that the  level set $L(\tau(x) )$ 
has positive Hausdorff dimension, i.e. 
$$\Gamma_H^{abs}:=\{ x \in [0,1]:~ \dim_{H}(L(\tau(x))) >0 \}.$$ 
Then $\Gamma_{H}^{abs}$ has 
full Hausdorff dimension, i.e.
\beql{183a}
\dim_{H}( \Gamma_H^{abs})= 1.
\eeq
\end{theorem}

We do not obtain information about the Lebesgue measure of $\Gamma_{H}^{abs}$.
Conjecture \ref{cg12} above is just the assertion  that $\Gamma_{H}^{abs}$
has Lebesgue measure $0$.\medskip

Theorem \ref{th13a}  is proved in Section \ref{sec6}. It is an immediate corollary of
an analogous result  for local level sets (Theorem \ref{th51}), which proves
Hausdorff dimension $1$ even when one restricts to $x \in \ddL$.
We define a family of  Cantor sets $\Lambda_{2r}$,
contained in $\Gamma_{H}^{abs} \cap \ddL$, 
where $r\ge 3$ is an integer parameter,
and show that  $\dim_{H}(\Lambda_{2r}) \to 1$ as $r \to \infty$.
Theorem \ref{th51}  also shows:
 {\em The deficient digit set $\ddL$ has Hausdorff dimension $1$.}\medskip

%
%
%
\begin{theorem}~\label{th13} {\em (Positive Hausdorff dimension ordinate level sets)}
Let $\Gamma_H^{ord}$ be  set of ordinates $y \in [0, \frac{2}{3}]$ such that the
Takagi function  level set $L(y)$ 
has positive Hausdorff dimension, i.e. 
$$\Gamma_H^{ord}:=\{ y:~ \dim_{H}(L(y)) >0 \}.$$ 
Then $\Gamma_{H}^{ord}$ has 
full Hausdorff dimension, i.e.
\beql{183}
\dim_{H}( \Gamma_H^{ord})= 1.
\eeq
\end{theorem}

Here we know that the  set $\Gamma_{H}^{ord}$ has 
 Lebesgue measure $0$, a result that   follows 
from Theorem~\ref{th12}. \\

Theorem \ref{th13} is proved in  Sect.  \ref{sec7}.
We show that the Takagi function 
$\tau(x)$ restricted to each Cantor  set $\Lambda_{2r}$  in Sect. \ref{sec6}
 is strictly increasing and
is a  bi-Lipschitz map. The result follows since 
bi-Lipschitz maps preserve Hausdorff dimension. \medskip

%
%

\subsection{Dimension Spectrum of Takagi Level Sets}\label{sec12}

Our results establish
 that level sets of the Takagi function
 exhibit a nontrivial dimension spectrum. 
We define the 
 {\em dimension spectrum function} for Takagi level sets to
be the function
\beql{169}
f_{\tau}(\alpha) := \dim_{H}\{ y:  \dim_{H} (L(y)) \ge \alpha\}.
\eeq
By definition this function is a nonincreasing function of $\alpha$.
The second author conjectured  \cite{Mad10}, and de Amo et al. \cite{ABDF11}
proved, that there are no level sets having Hausdorff dimension exceeding
$\frac{1}{2}$, which gives
\beql{170}
f_{\tau}(\alpha)=0 ~~~\mbox{for}~~\alpha > \frac{1}{2}.
\eeq
We trivially have $f_{\tau}(0)=1$ and Theorem~\ref{th13} of this paper establishes that
$\lim_{\alpha \to 0^{+}} f_{\tau}(\alpha) =1,$
by showing  that  for all sufficiently
large integers $r$,
\beql{171a}
f_{\tau}(\frac{1}{2r})  \ge 1- \frac{2\log r}{2r}.
\eeq
By monotonicity of this function  we have $f_{\tau}(\alpha)>0$ holding
on some  interval $[0, \alpha_0)$ of positive length.
It seems reasonable to expect that 
$f_{\tau}(\alpha) > 0 ~~\mbox{for}~~0 \le \alpha < \frac{1}{2}.$
Finally, we note that the  assertion that 
$f_{\tau}(\alpha) < 1$ holds for each $\alpha >0$ would imply the truth of 
Conjecture \ref{cg12}.

One can also study  analogous questions for  local level sets. 
We  define
the {\em  local level set dimension spectrum function} by
\beql{172a}
f_{\tau}^{\ast}(\alpha) := \dim_{H} \left( \{ x \in \ddL : \dim_{H} (L_x^{\loc}) \ge \alpha \} \right).
\eeq
This function is not directly comparable with $f_{\tau}(\alpha)$ because
it samples abscissa points rather than ordinate points. 
The nonexistence of level sets of Hausdorff
dimension exceeding $\frac{1}{2}$ yields
\beql{172b}
f_{\tau}^{\ast}(\alpha)=0~~~\mbox{for}~~\alpha > \frac{1}{2}.
\eeq
Theorem \ref{th51} of this paper establishes that 
$f_{\tau}^{\ast}(0) = \lim_{\alpha \to 0^{+}} f_{\tau}^{\ast}(\alpha) =1,$
by showing  that, for all
large enough  integers $r$, 
\beql{174}
f_{\tau}^{\ast}(\frac{1}{2r}) > 1- \frac{2\log r}{r}.
\eeq
Again it follows that  $f_{\tau}^{\ast}(x)>0$ on some interval $[0, \alpha_1)$
of positive length. \smallskip

The notion of {\em multifractal formalism} 
(or {\em thermodynamic formalism}) has been introduced
in connection with the H\"{o}lder  spectrum of points of
a nonsmooth function having given H\"{o}lder exponent;
see Jaffard \cite{Ja97a}, \cite{Ja97b} for discussion and many references.
It is predicted (under suitable hypotheses)
that the H\"{o}lder spectrum dimension $f(\alpha)$ of a given function
is a real-analytic function of $\alpha$ over a certain range,
and also exhibits a convexity property  over this range.
The multifractal formalism further relates local properties
of the function (local H\"{o}lder exponents) to global
smoothness properties of the function, invoving a Legendre transform
of the function $f(\alpha)$. 
Our dimension functions above encode
rather different properties of the Takagi function--size of its level sets--
but in a certain sense still sample local properties of this function.
Therefore it seems reasonable to
ask whether  the associated dimension functions $f_{\tau}(\alpha)$
and $f_{\tau}^{\ast}(\alpha)$ 
might have similar analytic and convexity properties.\smallskip

%
%

\subsection{Discussion} \label{sec13}

There has been much prior study of the non-differentiable structure
of the Takagi function and related functions according to various measures. 
See the papers of  
Allaart and Kawamura  (\cite{AK06}, \cite{AK10}),   and references therein.
We also mention the survey paper of the first author \cite{La11}.
Level sets study only one particular aspect of the nondifferentiability
of this function.\smallskip

An  interesting feature of the Takagi function
is that the cardinality  of ``generic" abscissa and ``generic" ordinate level sets
 in the Lebesgue measure sense differ drastically.
This difference can occur because
 sampling a point $x$ on the
abscissa favors  picking level sets which are ``large."  
It is a manifestation of the non-differentiable nature
of the Takagi function. 
This difference
 indicates that the Takagi function must 
 (in some sense) have ``infinite slope" over part of its domain.
 In particular $\tau(x)$ is not  a function of bounded variation.
 \smallskip
 
 We raise as open problems the determination of   
 the dimension functions $f_{\tau}(\alpha)$
and $f_{\tau}^{\ast}(\alpha)$ for $0 \le \alpha \le 1$ and of deciding
whether properties predicted by the thermodynamic formalism hold. 
It might conceivably be true that $f_{\tau}(\alpha) = f_{\tau}^{\ast}(\alpha)$
holds for $0 \le \alpha \le 1;$ the above discussion shows it holds
 for $\alpha \ge \frac{1}{2}.$ \smallskip

Some further improvements of the results on level sets 
appear in  Allaart \cite{A11}, \cite{A11b}. 
In particular \cite{A11b} obtains 
results on the detailed distribution of level sets having
a given finite number of elements.

\paragraph{Acknowledgments.}
The first author thanks D. E. Knuth for raising
  interesting questions about  the Takagi function
 (see \cite[Problem 82, p. 103]{Kn05}.)
We thank  Pieter
 Allaart for allowing use to include
 his substantially simplified proof of Theorem~\ref{th52}, 
 and for  bringing the work of
of Buczolich \cite{Buz08} to our attention.  We thank the reviewer
for helpful comments.

%
%
%
%
\section{Preliminaries: Properties of the Takagi Function}\label{sec2}
\setcounter{equation}{0}


\subsection{Functional Equations} \label{sec21}

We first recall
two  functional equations satisfied by the Takagi function \cite[Lemma 2.2]{LM10a}.
%
%
%
%

\begin{lemma}~\label{le22} {\em (Takagi functional equations)}
The Takagi function satisfies  two functional equations, valid for
$0 \le x \le 1,$ the reflection equation
\beql{206a}
\tau(x) = \tau(1-x),
\eeq
and the dyadic self-similarity equation
\beql{206b}
2 \tau(\frac{x}{2})= \tau(x) + x.
\eeq
\end{lemma}

We next formulate a local self -similarity property of the graph of the Takagi function.
To describe it  we  require some functions  determined by 
the binary expansion of $x$.

%
%
\begin{defi}\label{de23} 
{\em 
Let $x$ denote a binary expansion
\beql{221}
x := \sum_{j=1}^{\infty} \frac{b_j}{2^j}= 0.b_1 b_2 b_3..., 
\eeq
with each $b_j \in \{0, 1\}$. 
For each $j \ge 1$ we define the following integer-valued functions.

(1) The {\em digit sum function} $\dM_j(x)$ is
\beql{221a}
\dM_j(x) := b_1 +b_2 + \cdots + b_j.
\eeq
We also set
\beql{221b} 
N_j^{0}(x) :=  j - \dM_j(x).
\eeq
These functions count the number of $1$'s (resp. $0$'s) in the first
 $j$ binary digits of $x$. \smallskip

(2) The {\em  deficient digit function} $\dN_j(x)$ is given by
\beql{222}
\dN_j(x):=N_j^0(x) - \dM_j(x) = j- 2\dM_j(x)  = j-  2(b_1+b_2+ \cdots + b_j) .
\eeq
The function $\dN_j(x)$ counts the excess of binary digits $b_k=0$ over those with $b_k=1$
in the first $j$ digits, i.e. it is positive if there are more $0$'s than $1$'s. 
(Note that dyadic rationals have two different binary expansions, and the functions $N_j^0(x),
\dM_j(x)$, $\dN_j(x)$ depend
on the choice of expansion.)
}
\end{defi}

The   local
self-similarity property of the Takagi function graph is
as follows (\cite[Lemma 2.5(2)]{LM10a}).

%
%

\begin{lemma}\label{le25}
{\em (Takagi  self-similarity)} \\
 Let   $x_0 = 0.b_1b_2\ldots b_{n} 0^{\infty}= \frac {k}{2^{n}}$ 
 be a  dyadic rational, and parameterize the interval
 $[\frac{k}{2^n}, \frac{k+1}{2^n}]$ as 
 $x= x_0 + \frac{w}{2^n}$ for $0 \le w \le 1$.
Then there holds
\beql{251}
\tau(x) = \tau(x_0) + \frac{1}{2^n}\left(\tau(w) + D_n(x_0)\, w \right).
\eeq
That is, on the dyadic interval $ [\frac{k}{2^{n}}, \frac{k+1}{2^{n}}]$
the graph of the function $\tau(x)$ is a miniature version  of its full graph,
vertically shifted by $\tau(x_0)$, shrunk by a factor $\frac{1}{2^{n}}$,
and tilted by an additive linear factor $\frac{1}{2^n} D_n(x_0) w$.
\end{lemma}

In particular, for the case of balanced dyadic rationals,
which are ones with $D_{n}(x_0)=0$ (necessarily $n=2m$ is even),
\eqn{251} simplifies to
\beql{251b}
\tau( x) = \tau(x_0) + \frac{\tau(w)}{2^n},
\eeq
which comprises only a vertical shift and shrinking of the Takagi function.


\subsection{Local level sets} \label{sec22}

The notion of  local level set $L_x^{loc}$ is 
attached to the binary expansion of an abscissa point $x \in [0,1]$.
We show that certain combinatorial flipping operations applied
to the binary expansion of $x$ yield new points $x'$ in the same level set.
The totality of points reachable from $x$ by these combinatorial operations
will comprise the local level set $L_x^{loc}$ associated to $x$.

Let a binary expansion of $x \in [0,1]$ be given:
\beql{121}
x := \sum_{j=1}^{\infty} \frac{b_j}{2^j}= 0.b_1 b_2 b_3..., ~~~~ \mbox{each}~b_j \in \{0, 1\}. 
\eeq
The {\em flip operation} (or {\em complementing operation}) on  a single binary digit $b$ is 
\beql{124a}
\bar{b} := 1- b.
\eeq
We associate  to any binary expansion $x$  the sequence  of  digit positions $j$
at which tie-values of the deficient
digit function $\dN_j(x)=0$ occur, which we call {\em balance points}; note that all such
$j$ are even.  The  {\em balance-set}  $Z(x)$ associated to $x$ is denoted 
\beql{123}
Z(x) : = \{ c_k:~~\dN_{c_k}(x)=0\}.
\eeq
where we define $c_0=c_0(x) = 0$ and set $c_0(x)< c_1(x)< c_2(x) < ...$. This sequence of tie-values may
be finite or infinite. If it is finite, ending in $c_{n}(x)$, we make the
convention to adjoin a final ``balance point"  $c_{n+1}(x)= +\infty$.
 We call a {\em ``block"}  an indexed set of digits between two consecutive balance points, 
\beql{124}
B_k(x) := \{ b_j: ~c_k(x)  < j \le c_{k+1}(x)\},
\eeq
which includes the second balance point but not the first.
We define an equivalence relation on
blocks, written  $B_k(x) \sim B_{k'}(x')$ to mean the block endpoints agree
($c_k(x)= c_{k'}(x')$ and $c_{k+1}(x) = c_{k'+1}(x')$)
and either $B_k(x) = B_{k'}(x')$ or $B_k(x) = \bar{B}_{k'}(x')$, where the bar operation flips
all the digits in the block, i.e. 
\beql{124c}
b_j \mapsto \bar{b}_j:= 1- b_j,~~~~~~~ c_k < j \le c_{k+1}.
\eeq
Finally, we define the equivalence relation $x \sim x'$
on two binary expansions  to mean that
they have identical balance-sets 
$Z(x) \equiv Z(x')$, and furthermore every block $B_k(x) \sim B_k(x')$ for $k \ge 0$.
Note that $x \sim 1-x$; this corresponds to a flipping operation being applied to
every binary digit. 
In \cite{LM10a} we
showed that the  equivalence relation $x\sim x'$ implies that $\tau(x)= \tau(x')$,
 so that $x$ and $x'$
are in the same level set of the Takagi function.

%
%

\begin{defi}~\label{de12}
{\em 
The {\em local level set} $L_x^{loc}$ associated to $x$ is the set of equivalent points, 
\beql{125}
L_x^{\loc} := \{ x': ~~x' \sim x\}.
\eeq
We use again the convention that $x$ and $x'$ denote binary expansions, so
that dyadic rationals require special treatment. }
\end{defi}

We recall some basic properties of  local level sets 
(\cite[Theorem 3.1, Corollary 3.2]{LM10a}).

\begin{enumerate} 
\item
Local level sets  $L_x^{\loc}$ are closed sets. Two local level sets
either coincide or are  disjoint. 

\item
  Each local level set $L_x^{\loc}$ is contained in a level set:
 $L_x^{loc} \subseteq L(\tau(x))$.
That is, if  $x_1 \sim x_2$ then $\tau(x_1) = \tau(x_2).$

\item
 Each level set $L(y)$ partitions into local level sets
\beql{300e}
L(y) = \bigcup_{{x \in \ddL}\atop{\tau(x)= y}}  L_{x}^{\loc}
\eeq
Here $\ddL$ denotes the collection of leftmost endpoints of all local level sets.

\item
 A local level set $L_x^{loc}$ is a  finite set if the balance-set $Z(x)$ is finite;
otherwise it is an uncountable perfect set (Cantor set).
\end{enumerate}


\subsection{Deficient digit set $\ddL$} \label{sec23}

In \cite{LM10a} 
we  studied the set of leftmost endpoints $\ddL$ of local level sets.
 \smallskip

%
%

\begin{defi}\label{de31}
{\em
 The {\em deficient  digit  set} $\ddL$ consists of all
$x$ such that
$$
\ddL := \{ x = \sum_{j=1}^{\infty} \frac{b_j}{2^j}:  ~~ \dN_j(x) \ge 0 ~\mbox{for ~all} ~~j \ge 1\}.
$$
}
\end{defi}

The deficient digit set is  a Cantor-type set obtained by removing a
certain countable collection of open intervals from the unit interval, which we describe
using the following definitions.

%
\begin{defi} \label{de32}
{\em
(1) The {\em breakpoint set}
 $\sBp$ consists  of $\sBo'=0$ together with the  collection of  all balanced
 dyadic rationals in $\ddL$. These are all 
   $B' = \frac{n}{2^{2m}}$  that have  binary expansions of the form
$$
B'=0.b_1 b_2 ... b_{2m-1}b_{2m} \,  0^{\infty}~~~\mbox{for some}~~ m \ge 1,
$$
that satisfy  the condition
\beql{331b}
\dN_j(B') \ge 0 ~~~\mbox{for} ~~1 \le j \le 2m-1, ~~~\mbox{and}~~~D_{2m}(B') =0,
\eeq
This condition implies $b_{2m}=1$.\\

(2) The {\em small  breakpoint set } $\sB$ is the subset of 
the breakpoint set $\sB'$ consisting of
$\sBo=0$ plus all dyadic rationals in $\sB'$ that satisfy the extra condition that 
the last two binary digits $b_{2m-1}= b_{2m}=1.$

}

\end{defi}

In \cite{LM10a} we used the small breakpoint set $\sB$
to label  the intervals removed from
$[0,1]$ to create the deficient digit set $\ddL$. 

%

\begin{defi} \label{de33b}
{\em
For each dyadic rational  $B=0.b_1b_2 ...b_{l} 01^k$, $k \ge 2$ 
 in the small breakpoint set $\sB$ ($B \ne \sBo$) associate
 the open interval
$I_B := (x( B)^{-}, x(B)^{+})$
 having endpoints
\begin{eqnarray*}
x(B)^{-} & := &0. b_1 b_2 ... b_l 0 1^k (01)^{\infty}\\
x(B)^{+} & := & 0. b_1 b_2 ... b_{l} 1 0^k (00)^{\infty},
\end{eqnarray*}
necessarily with $k \ge 2$. We also set $I_{\sBo} = (\frac{1}{3}, 1).$
} 
\end{defi}

The following result gives properties of the deficient digit set $\ddL$
( \cite[Theorem 4.6]{LM10a}.)

%
%

\begin{theorem}\label{th29} {\em(Properties of  Deficient Digit Set)}\\
$~~~~~$(1) The deficient digit  set $\ddL$ comprises the set of leftmost endpoints of all local
level sets. It satisfies $\ddL \subset [0, \frac{1}{3}]$.

\smallskip
(2) The deficient digit sum set $\ddL$ is a closed, perfect set (Cantor set).
It is given by 
\beql{342}
\ddL= [0, 1) \backslash \bigcup_{ B \in \sB} I_{B}.
\eeq
where the omitted open intervals $I_B$, for $B$ in
the small breakpoint set, have  right endpoint  a dyadic rational and 
left endpoint a rational number with denominator $3 \cdot 2^k$ for some $k \ge 1$.

\smallskip
(3) The deficient digit set $\ddL$ has Lebesgue measure zero.
\end{theorem}

In \cite[Lemma 4.5]{LM10a} it is shown that the value of the
endpoints of the removed intervals satisfies
\beql{271}
x_B^{+} - x_{B}^{-}  = \tau(x(B)^{-})- \tau(x(B)^{+})= \frac{1}{2^{k+l} \cdot 3},
\eeq
so that linear interpolation of a function across the a removed interval
always has slope $-1$.

%
%
%

\subsection{Takagi function on deficient digit set}\label{sec32}

In \cite[Theorem 4.8]{LM10a} we proved that the
Takagi function is nondecreasing on the  set $\frac{1}{2}\ddL$. Note
that
$$
\frac{1}{2} \ddL:= \{ \frac{1}{2} x: x \in \ddL\}
 =\{ x \in
[0,1]: \dN_j(x) > 0 \textrm{ for all } j \geq 1\},
$$
which shows that $\frac{1}{2} \ddL \subset \ddL$.
%
%

\begin{theorem}\label{th33} 
  (1) The Takagi function  is nondecreasing on the set $\frac{1}{2} \ddL$.
  
   (2) The Takagi function  is strictly increasing on $\frac{1}{2} \ddL$
  away from a countable set of points, which are a subset of  those rationals having binary
  expansions ending in $0^{\infty}$ or $(01)^{\infty}$. For each level
  $y$ the  equation $y= \tau(x)$
   has at most two solutions with $x \in \frac{1}{2} \ddL$. Thus  if $x_1< x_2< x_3 $ all
  belong to  $\frac{1}{2}\ddL$ then $\tau(x_3) > \tau(x_1)$.
   \end{theorem}

This result will be used in establishing 
the bi-Lipschitz property appearing in Theorem \ref{th52}.

%
%
%
%

\section{Takagi Singular Measure} \label{sec3}

A main ingredient used in this paper will be a singular measure,
which is the weak derivative of
 the Takagi singular function constructed in \cite{LM10a}.
 We summarize the basic facts on the Takagi singular function,
 taken from  \cite[Theorems 1.5 ]{LM10a}. It is pictured in
 Figure \ref{fig421} below.

%
%
%
\begin{theorem} \label{th15a} {\em (Takagi singular function)}
The function  $\tauS(x)$ defined by $\tauS(x)= \tau(x) + x$ for
 $x \in \ddL$  
is a nondecreasing function on $\ddL$. Define its extension to
all $x \in [0,1]$ by 
$$
\tauS(x) := \sup\{ \tauS(x_1):  x_1 \le x ~~~\mbox{with}  ~~ x_1 \in \ddL\}.
$$
 Then the  function $\tauS(x)$  is  a monotone singular function. That is, it is
 a nondecreasing continuous function  
 having  
 $\tauS(0)=0, \tauS(1)=1$, which has derivative zero at (Lebesgue) almost all points of $[0,1]$. 
  The closure of the set of points of increase of $\tauS(x)$ is the deficient digit set $\ddL$.
\end{theorem}

\smallskip

%
%

\begin{figure}[h]
  \begin{center}$                                                               
    \begin{array}{cc}
      \includegraphics[height=2in]{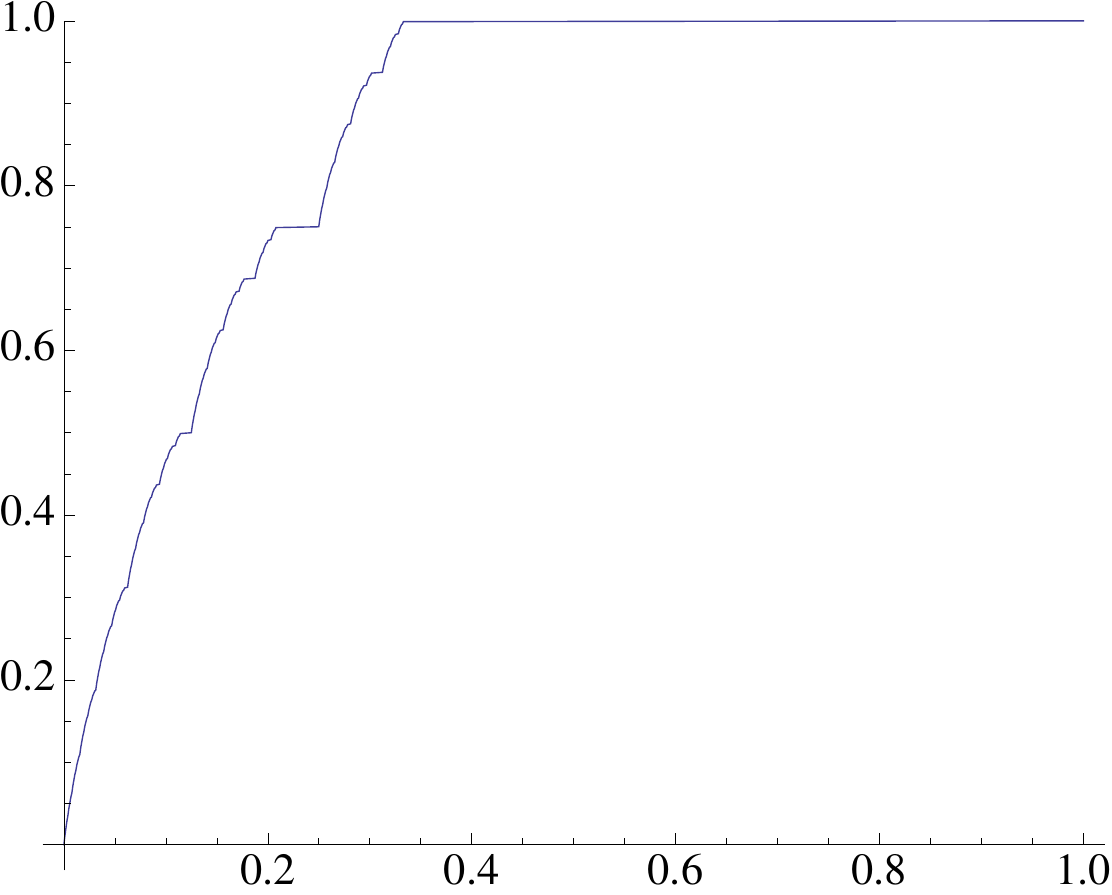}
    \end{array}$
  \end{center}
\caption{Graph of 
 the
  singular Takagi function.}
\label{fig421}
\end{figure}

The weak derivative of
the Takagi singular function is a Radon measure
$\muT$ that we call the {\em Takagi singular measure},
and it has the following properties.

%
%
\begin{theorem}\label{th34} {\em (Takagi singular measure)}
 The  Takagi singular function $\tauS(x)$  is the 
 definite integral 
 of a  nonnegative Radon 
 measure $\muT$, i.e.
$\tauS(x) = \int_{0}^x d\muT.$
The measure $\muT$ has support equal to the deficient digit set $\ddL$, so that
\beql{332b}
\int_{0}^1 d\muT= \int_{\ddL} d\muT = 1.
\eeq
Thus it  is singular continuous with respect to Lebesgue measure.
For every Borel set $K$ in $[0,1]$, 
\beql{332c}
\muT(K) = \meas( \tauS(K)),
\eeq
where $\meas(\cdot)$ denotes Lebesgue measure.\smallskip
\end{theorem}

%
%

\paragraph{Proof.}
By  Theorem~\ref{th15a},
 the  function $\tauS(x)$ is monotone and bounded, hence of 
 bounded (pointwise) variation.
 It follows   that 
 its distributional derivative is a finite Radon measure, call it $\muT$. It is
 necessarily nonnegative
 since $\tauS(x)$ is nondecreasing (cf. \cite[Theorem 1, Sect. 5.1]{EG92}).
 
 The support of a measure $\mu$ on 
 the real line is the closure of the set of points of 
 increase of $\mu$, which is $\ddL$ by Theorem  \ref{th15a}.
 Since $\meas(\ddL)=0$ it is a singular measure, and
  it is singular continuous because the
 integral $\tauS(x) := \int_{0}^{x} d\muT$ is a continuous function of $x$.
  Finally, since  $[0,1]$ is compact, any Radon measure (inner regular
 measure) on $[0,1]$  is also outer regular
 (\cite[Theorem E, Sect. 52]{Hal50}).
 Outer regularity means that for any Borel set $K$,
 $\muT(K) = \inf_{ K \subset U} \{\muT(U): ~U ~~\mbox{an open set}\}.$
 We know that the Takagi singular measure of an open interval $U= (0, t)$ is
 $$
 \muT(U) := \int_{0}^t d \muT(x)=  \tauS(t) - \tauS(0)= \tauS(t) = \meas( (\tauS(0), \tauS(t))) = \meas( \tauS(U)).
 $$
 where $\meas$ denotes Lebesgue measure.
 Outer regularity of both $\muT$ and of Lebesgue measure now implies that  
 the equality holds on every Borel set $K$, proving \eqn{332c}. $~~~\Box$\\

%
%
%
%
\section{Structure of Takagi Singular Measure}\label{sec4}
\setcounter{equation}{0}

We now compute values of the Takagi singular measure
on various subsets of $\ddL$ In \S\ref{sec24} we define a fine partition of $\ddL$
and in the remaining subsections 
we compute the singular measures of all sets in this partition.

\subsection{Fine partition of deficient digit  set} \label{sec24}

We  partition the set $\ddL$ of left endpoints of local level sets into finer pieces, as
follows:
\beql{401}
\ddL= \ddL_{\infty} \, \bigcup \,  \ddL_{fin}
\eeq
in which 
\beql{401d}
\ddL_{\infty} := \{ x \in \ddL:  ~\dN_j(x) = 0 ~\mbox{for infinitely many}~ j \ge 1\}.
\eeq
and $\ddL_{fin}$ is its complement, 
\beql{401e}
\ddL_{fin} := \{ x \in \ddL:  ~\dN_j(x) \ge 1 ~\mbox{for all sufficiently large}~ j \}.
\eeq
The latter set can be further partitioned into subsets labelled by
elements $B' $ in the breakpoint set  $\sB'$
consisting of all dyadic rationals $B'= 0.b_1b_2 ..b_{2m} = \frac{k}{2^{2m}}$ such
that all $\dN_j(B') \ge 0$ and $\dN_{2m}(B') =0$.
To each $B'\in \sB'$ we associate the set
\beql{403a}
\ddL(B'):= \{ x= B' + \frac{x'}{2^{2m}}: ~x' \in \frac{1}{2} \ddL\}.
\eeq
In particular  for $m=0$ we have one set $B'=B_0=0$ with  $\ddL(B_0) = \frac{1}{2} \ddL$.

%
%

\begin{lemma}~\label{le41} {\em (Fine Partition of Deficient Digit Set)}
 The set $\ddL_{fin}$ has a partition
\beql{403b}
\ddL_{fin} = \bigcup_{ B' \in \sBp} \ddL(B'),
\eeq
with union over the breakpoint set $\sBp$.
Each set $\ddL(B')$  is a closed set. 
\end{lemma}

\paragraph{Proof.} 
Elements  $x \in \ddL(B')$  have $\dN_j(x) \ge 0$ for all $j \ge 1$,
$\dN_{2m}(x)=\dN_{2m}(B')=0$, and $\dN_{j}(x)>0$ for $j \ge 2m+1$. In particular
$\ddL(B') \subset \ddL_{fin}$. The sets are disjoint for different $B'$
because the value $2m$ is uniquely determined for each element of the
set $\ddL(B')$, and this determines the initial digits $B'$ uniquely. 
Finally we see that each element $x \in \ddL_{fin}$ has 
 associated to it  a unique maximal value $2m$  of $j$ such that $\dN_j(x)=0$,
($j$ is necessarily even) and this assigns it to a particular $\ddL(B')$. 
$~~~\Box$\\


\subsection{Singular measure mass calculations} \label{sec41}

The Takagi singular  measure $d\muT$ is not translation-invariant.
However in  Theorem \ref{th41} below we are able to use  its self-similarity
properties to compute the $\muT$-measure of certain sets inside $\ddL$,
namely the sets $\ddL_{\infty}$ and the sets $\ddL(B')$ in the fine partition
of $\ddL$ in Lemma \ref{le41}, where $B'$ runs over 
the set of balanced dyadic rationals
that belong to $\ddL$.

%
%

\begin{theorem}~\label{th41}
 For each balanced dyadic rational $B' =0.b_1 b_2... b_{2m}= \frac{k}{2^{2m}}$ in
 the deficient digit set  $\ddL$ 
 the fine partition set  $\ddL(B')$ is a closed set, and its  Takagi singular measure is
\beql{411} 
\muT( \ddL(B')) := \int_{\ddL(B')} d\muT = \frac{1}{2^{2m+1}}.
\eeq
\end{theorem}

\paragraph{Proof.} 
We already know that  $\muT(\ddL) = 1$ via Theorem~\ref{th34}.

\paragraph{Claim 1.} {\em The Takagi singular measure mass of $\frac{1}{2} \ddL$ is
given by
\beql{419}
\muT( \ddL(B_0)) := \muT( \frac{1}{2} \ddL) = \frac{1}{2} \muT(\ddL) = \frac{1}{2}.
\eeq
}

To prove the claim, we use the self-similarity relation
 in Lemma \ref{le22}(1), 
$$
2\tau(\frac{1}{2} x) = x + \tau(x),~~~~~\mbox{for} ~~~0 \le x \le 1.
$$
It $x \in \ddL$ then $\tau(x) = \tauL(x)$ so that we obtain
\beql{426a}
2 \tau(\frac{1}{2} x) = x + \tau(x) 
= \tauS(x).
\eeq
Thus if $x_1 <x_2$ with both $x_i \in \ddL$, then
\begin{eqnarray}\label{427a}
\int_{\frac{1}{2} x_1}^{\frac{1}{2} x_2} \muT &=& \tauS(\frac{1}{2} x_2) - \tauS(\frac{1}{2} x_1) \nonumber\\
&=& \left( \tau(\frac{1}{2} x_2) + \frac{1}{2} x_2\right) -
\left( \tau(\frac{1}{2} x_1) + \frac{1}{2} x_1\right) \nonumber \\
&=& \frac{1}{2}\left( \tauS(x_2) - \tauS(x_1) \right) + \frac{1}{2} \left( x_2 - x_1\right)  \nonumber \\
&=& \frac{1}{2} \int_{x_1}^{x_2} \muT +   \frac{1}{2} \left( x_2 - x_1\right).
\end{eqnarray}
We may rewrite this as
\beql{428a}
\left|\int_{\frac{1}{2} x_1}^{\frac{1}{2} x_2} \muT -  \frac{1}{2} \int_{x_1}^{x_2} \muT \right| 
\le \frac{1}{2} \meas([x_1, x_2])
\eeq
where the last term denotes the Lebesgue measure of the interval $[x_1, x_2]$. 

Now by \eqn{342} for each $m \ge 1$ we obtain a covering of $\ddL$ using
\beql{428b}
\ddL \subset \sP_{2m} := [0,1) \smallsetminus \bigcup_{{B \in \sB}\atop{|B| \le 2m}} I_{B}.
\eeq
in which we remove only a finite number of the  ``flattened"  open intervals $I_{B}$ corresponding
to those $B \in \sB$ (the small breakpoint set) having dyadic length at most $2m$. The set $\sP_{2m}$ is
a closed set comprised of a finite number of intervals, $[x_j, x_j^{'}]$, say, 
having both endpoints $x_j, x_j^{'} \in \ddL$. Adding up the relations \eqn{428a} over
these intervals yields
\beql{429}
\left|\int_{\frac{1}{2} \sP_{2m}} \muT -  \frac{1}{2} \int_{\sP_{2m}} \muT\right| \le \frac{1}{2} \meas(\sP_{2m}),
\eeq
in which $\meas(\sP_{2m})$ denotes the Lebesgue measure of the set $\sP_{2m}$. 
Next we note that the $\sP_{2m}$ form 
a nested family  $\sP_2 \supset \sP_4 \supset \sP_{6} \supset \cdots$ of closed sets,
with
$$
\ddL = \bigcap_{m=1}^{\infty} \sP_{2m}.
$$
Since these sets are Borel measurable 
and  $\muT$  is outer regular 
we have 
$$
\lim_{m \to \infty} \int_{\sP_{2m}} \muT = \int_{\ddL} \muT,
$$
$$
\lim_{m \to \infty} \int_{\frac{1}{2}\sP_{2m}} \muT = \int_{\frac{1}{2}\ddL} \muT,
$$
cf. Evans and Gariepy \cite[Theorem 1, p. 2]{EG92}.
Now  Theorem~\ref{th29} (2), (3) 
 together establish that
$$
\meas(\sP_{2m}) \to 0 ~~~\mbox{as}~~~m \to \infty.
$$
Thus letting $m \to \infty$ in \eqn{429} yields
$$
\int_{\frac{1}{2}\ddL} \muT= \frac{1}{2} \int_{\ddL} \muT,
$$
which with $\int_{\ddL} \muT =1$ proves Claim 1.

\paragraph{Claim 2.} {\em Let $B' = \frac{k}{2^{2m}} \in \sBp$ and
suppose that $x_i = B'+ \frac{x_i'}{2^{2m}}$ for $i=1,2$, with both  $x_i^{'} \in \ddL$. 
Then
\beql{432a}
 \int_{x_1}^{x_2} \muT = \frac{1}{2^{2m}} \left(\int_{x_1^{'}}^{x_2^{'}} \muT  + (x_2^{'} - x_1^{'})\right).
\eeq
}

Since $B'$ is a balanced dyadic
rational, we may deduce  \eqn{432a} using the formula of Lemma~\ref{le25}, in analogous fashion
to \eqn{427a}. This proves Claim 2. \smallskip

Now we complete the proof. 
Claim 2 yields the formula
\beql{433a}
\left| 2^{2m} \int_{x_1}^{x_2} \muT - \int_{x_1^{'}}^{x_2^{'}} \muT \right| \le \meas ([ x_1^{'}, x_2^{'}]).
\eeq
For any $B' \in \sB'$ we have
$$
\ddL(B') := \{ x:~x=B' + \frac{x'}{2^{2m}} ~~\mbox{with}~~x^{'} \in \frac{1}{2} \ddL\}.
$$
Now we may cover  the set $\frac{1}{2} \ddL$ with $\frac{1}{2} \sP_{2n}$.
We apply 
the approximation bound \eqn{433a} summed up over all the intervals in
$\sP_{2n}$, to obtain
$$
\left|2^{2m} \int_{ B' + \frac{1}{2^{2m}} (\frac{1}{2} \sP_{2n})} \muT -
\int_{\frac{1}{2}\sP_{2n}} \muT \right| \le \frac{1}{2} \meas(\sP_{2n}).
$$
Letting $n \to \infty$ we deduce,
 using $\meas(\sP_{2n}) \to 0$, that
$$
2^{2m}\int_{ B' + \frac{1}{2^{2m}} (\frac{1}{2} \ddL)} \muT=
\int_{\frac{1}{2}\ddL} \muT = \frac{1}{2}.
$$
This yields, since 
$\ddL(B') := B'+ \frac{1}{2^{2m}}( \frac{1}{2}\ddL) $,  that
$$
\int_{\ddL(B')} \muT= \frac{1}{2^{2m+1}},
$$
as asserted. $~~~\Box$.\\


\subsection{Singular measure mass of $\ddL_{\infty}$} \label{sec61}

The calculations of the last section yield the following result.

%
%

\begin{theorem}~\label{th42} {\em (Takagi singular measure: fine partition)}
Let $\muT$ denote the Takagi singular measure. 
The sets $\ddL_{fin}$ and $\ddL_{\infty}$ are Borel sets,
hence measurable. 
We have 
\beql{421}
\muT(\ddL_{fin}) = 1, 
\eeq
which shows that
\beql{422}
\muT(\ddL_{\infty}) = 0.
\eeq
Consequently the  image of this set under the Takagi singular function $\tauS$ satisfies
\beql{422a}
\meas ( \tauS( \ddL_{\infty})) = 0,
\eeq
where $\meas$ denotes Lebesgue measure.
\end{theorem}

\paragraph{Proof.}
Each set $\ddL(B')$ is closed, hence their disjoint union 
$\ddL_{fin}$ is a Borel set, hence is $\muT$-measurable. The set $\ddL$ is closed, hence
$\ddL_{\infty} = \ddL \smallsetminus \ddL_{fin}$ is also 
a Borel set, hence is $\muT$-measurable, and
$$
\muT( \ddL_{\infty} )= \muT(\ddL) - \muT(\ddL_{fin}).
$$
(In fact one can easily show  that 
the closure of $\ddL_{fin}$ is $\ddL$.)

Since $\muT(\ddL)=1$ the theorem will follow on  showing $\muT(\ddL_{fin})=1$.
We have
$$
\muT ( \ddL_{fin} ) = \sum_{ B  \in \sBp}  \muT (\ddL(B)),
$$
where $\sBp$ is the breakpoint set. 
Theorem~\ref{th41} now gives $\muT(\ddL(B)) = \frac{1}{2^{2m+1}}$,
where $B = 0.b_1 \cdots b_{2m} = \frac{k}{2^{2m}}$, with $k$ odd. 
Recall from \cite[Lemma 4.2]{LM10a} that the number of balanced dyadic
rationals in $\ddL$ having the form
$\frac{k}{2^{2m}}$ for an odd $k$ is  the $m$-th Catalan number
$C_m = \frac{1}{m} \left( {{2m}\atop{m}} \right)$. Here for $m=0$ we
have $C_0=1$ corresponding to the element $B_0=0.$ 

The Catalan numbers are well known to have the generating function
\beql{330b}
\sum_{j=0}^{\infty} C_m z^{2m} = \frac{1- \sqrt{1-4z^2}}{ 2z^2}.
\eeq 
In consequence, taking $z= \frac{1}{2}$, we obtain
$\sum_{j=0}^{\infty} \frac{C_m}{2^{2m}} = 2.$
Therefore we obtain, using Theorem \ref{th41}, that 
$$
\muT( \ddL_{fin} ) = \sum_{m=0}^{\infty} C_m \frac{1}{2^{2m+1}}= \frac{1}{2}\left( \sum_{m=0}^{\infty}
 \frac{C_m}{2^{2m}}\right)= 1,
$$
which proves \eqn{421}. Now \eqn{422} follows,
and \eqn{422a} follows from Theorem \ref{th34} on taking $K= \ddL_{\infty}$
in \eqn{332c}.  $~~~\Box$\\

%
%
%
%
\section{Cardinality of Global Level Sets}\label{sec5}
\setcounter{equation}{0}

In this section we prove  Theorem~\ref{th12}, which states that for a full measure set
of ordinates $y$ the level set $L(y)$ is a finite set, and  that the
expected number of elements in  this set, with respect to Lebesgue measure on $0 \le y \le \frac{2}{3}$,
 is infinite. 
 
 We use the following result \cite[Theorem 5.8]{LM10a} giving
 the expected number of local level sets in a uniformly
 chosen level $y$ in $[0, \frac{2}{3}]$.

%
%
\begin{theorem}\label{th37} {\em (Expected number of local level sets)}
For a full Lebesgue measure set of ordinate points $y \in [0, \frac{2}{3}]$ the
number $N^{loc}(y)$ of local level sets at level $y$ is finite. Furthermore
\beql{361}
\int_{0}^{\frac{2}{3}}  N^{loc}(y) dy = 1.
\eeq
That is,  the expected number of local level sets  on 
 a randomly drawn ordinate level $y$  is $\frac{3}{2}$.
\end{theorem}

\noindent Theorem ~\ref{th37} was proved using  the coarea 
formula for functions of bounded variation
applied   to the flattened Takagi function $\tauL$.

%
%
\paragraph{Proof of Theorem~\ref{th12}.}
(1) Let $\Gamma_{\infty}^{ord}$  be the set of infinite levels, i.e. 
\beql{500a}
\Gamma_{\infty}^{ord} := \{ y: ~~L(y) ~~\mbox{ is an infinite set} \}.
\eeq
To show a full measure set of ordinates have finite level sets,
we show the contrapositive, that $\Gamma_{\infty}^{ord}$ has
Lebesgue measure $0$.  We have
$$
\Gamma_{\infty}^{ord} \subset \tau(\ddL_{\infty}) \, \bigcup \, \Lambda_{\infty}^{loc},
$$
in which 
$\tau(\ddL_{\infty}):= \{ y= \tau(x): ~ x \in \ddL_{\infty} \}$
detects all levels that contain at least one  infinite local level set,
and 
\beql{500f}
\Lambda_{\infty}^{loc} := \{ y: L(y)~~\mbox{contains infinitely many different local level sets}\}.
\eeq
 Now \cite[Theorem 7.2 (1)]{LM10a} shows that
$\Lambda_{\infty}^{loc}$ has Lebesgue measure $0$. Thus it suffices
to prove that $\tau(\ddL_{\infty})$ has Lebesgue measure $0$.

Recall that 
$
\tauS(x) = \tau(x) +x, ~~\mbox{for}~~x \in \ddL.
$
 Now consider $\tauS(x)$ restricted to $x \in \ddL(B)$
for a particular $B \in \sBp$, the   breakpoint   set.
We write $B= 0. b_1 b_2 \cdots b_{2m} =  \frac{k}{2^{2m}} $ where
$k$ is necessarily odd. Then $x \in \ddL(B)$ if and only if 
$$
x = B + \frac{\frac{1}{2}x'}{2^{2m}},~~~~\mbox{with}~ \frac{1}{2}x' \in \frac{1}{2} \ddL.
$$
Lemma \ref{le25} then gives
$$
\tau(x) = \tau(B) + \frac{1}{2^{2m}} \tau(\frac{1}{2}x'), ~~\mbox{with} ~~x' \in \ddL.
$$
We recall from Lemma \ref{le22} that $2\tau(\frac{1}{2}(x)) = \tau(x) + x$ if $x \in \ddL$,
whence
\beql{eq713}
2^{2m+1}\left(  \tau(x) - \tau(B)\right)= \tau(x')+ x' 
= \tauS(x'), ~~~\mbox{for}~ x' \in \ddL.
\eeq
Now  the linear map
$$
y \mapsto y' :=2^{2m+1}\left(y- \tau(B) \right)
$$
sends the interval $[\tau(B), \tau(B) + \frac{1}{2^{2m+1}}]$ onto $[0, 1]$ and 
it follows from the above that it 
sends $\tau(\ddL(B)) \subset [\tau(B), \tau(B) + \frac{1}{2^{2m+1}}] $
 onto  the range
$\tauS(\ddL)=[0,1]$. 
We see from the linearity of this map that
$
\tau(\ddL(B) )= [\tau(B), \tau(B) + \frac{1}{2^{2m+1}}],
$
$$
\meas \left( \tau(\ddL(B))\right) = 
\frac{1}{2^{2m+1}}.
$$
Adding up these contributions, the summation in  
Theorem \ref{th42} gives that 
 the total Lebesgue measure in $y \in [0, \frac{2}{3}]$
covered by images of these sets (counting overlaps with multiplicity)  is 
$$
\sum_{B \in \sBp} \meas(\tau( \ddL(B))) = \sum_{B \in \sBp} \muT(\ddL(B)) = 1.
$$
(The images have some overlap, allowing their total measure to exceed the length
of the interval $[0, \frac{2}{3}]$.) 
Viewing these points $x  \in \ddL(B)$ as 
labelling left endpoints of local
level sets, this says that a 
lower bound
 of the total number of local level set endpoints Lebesgue-integrated over $0 \le y \le \frac{2}{3}$,
 counted with multiplicity,  is $1$. Here  we did not count any local level set endpoints in
$\tau (\ddL_{\infty} ) := \tau\left( \ddL \smallsetminus \ddL_{fin}\right)$, coming
from the image of $\ddL_{\infty}$, since each $\ddL(B)$ is disjoint from $\ddL_{\infty}$. 
Theorem \ref{th37} now gives
$$
\int_{0}^{\frac{2}{3}} N^{\loc}(y) dy =1,
$$
where  $N^{loc}(y)$  counts the number of all local level set endpoints. 
We have already accounted for the full mass of this integral above, and any
omitted contribution to $N^{loc} (y)$ coming from 
$\tau (\ddL_{\infty} ) := \tau\left( \ddL \smallsetminus \ddL_{fin}\right)$ 
necessarily contributes an additional nonnegative amount.
Thus we may conclude that
$$
\meas( \tau(\ddL_{\infty})) = 0.
$$
as asserted.  \smallskip

(2) We aim to show that the expected size of  a global  level set is infinite, i.e. to show that
$$
\int_{0}^{\frac{2}{3}} | L(y)| dy = +\infty,
$$
where $|L(y)|$ counts the number of elements in $L(y)$. 
By the discussion above we have
\beql{451}
\int_{0}^{\frac{2}{3}} | L(y)| dy =\int_{0}^1 |L(\tau(x))| \muT(x) \geq
\sum_{B \in \sB'} \frac{1}{2^{2m+1}} 2^{r(B)}, 
\eeq
in which 
$$
r(B):= | \{ 1\le j < \infty: ~ N_j(B)=0\}|.
$$
Here each $r(B)$ is finite and bounded above by $m$ if $B= 0. b_1 b_2 \cdots b_{2m}$.
We rewrite this as
\beql{451b}
\int_{0}^{\frac{2}{3}} | L(y)| dy =\sum_{m=0}^{\infty}  \frac{L_m}{2^{2m+1}},
\eeq
in which
$$
L_m := \sum_{B \in \sB', |B| = 2m}   2^{r(B)}.
$$
Now we observe
that $L_m$, the total number of binary sequences of length $2m$ having $N_{2m}(B)=0$,
 has a combinatorial interpretation
as counting the number  the two-dimensional lattice paths of length $2m$ 
starting at the origin $(0,0)$,  taking steps either $(1,1)$ or $(1,-1)$,
and ending at $(2m, 0)$. Indeed these paths  groups into collections of paths of
size $2^{r}$ under the ``flipping" (reflection) operation, with each group containing exactly
one path in $\sB'$. (See the discussion and proof in Feller \cite[Theorem 4, p. 90]{Fe68}
and  also \cite[Lemma 4.2]{LM10a}.)
It follows that 
$$
L_m = \left( {{2m} \atop{m}}\right) \sim \frac{1}{\sqrt{\pi m}} 2^{2m}.
$$
Thus the terms in the series on the right side of \eqn{451b}
decay like $\Omega ( \frac{1}{\sqrt{m}})$, so the series \eqn{451b}
diverges, giving the result. $~~~\Box$. \\

%
%
%
%

\section{Level Sets of  Positive Hausdorff Dimension: Abscissa View}\label{sec6}
\setcounter{equation}{0}

We  study level sets having positive Hausdorff dimension.
In  the paper 
\cite[Sect. 3.3]{LM10a}) we classified those local level sets containing a rational number $x$
that are of positive Hausdorff dimension: this gives some explicitly determinable 
rational ordinates $y$ having this property. Here we show that  the set $\Gamma_H^L$ of
abscissa points in $\ddL$ that give local level sets having positive Hausdorff
dimension has full Hausdorff dimension $1$.  

%
%

\begin{theorem}\label{th51} {\em (Local level sets of positive Hausdorff dimension)}
Let $\Lambda_H^L$ denote the  set of $x \in \ddL$ such that the Hausdorff dimension 
of $L_x^{\loc}$
is positive, i.e.
\beql{600a}
\Lambda_{H}^L:= \{ x \in \ddL: ~~dim_H(L_x^{loc}) > 0\}.
\eeq
This set has full Hausdorff dimension, i.e. 
\beql{600b}
dim_H ( \Lambda_H^L ) = 1.
\eeq
 In particular, the deficient
digit set $\ddL$ has Hausdorff dimension $1$.
\end{theorem}

\paragraph{Proof.}  It clearly suffices to prove the first assertion.
For integer $r \ge 1$ let $\Gamma_{2r} $ consist of all abscissas $x \in [0,1]$
that satisfy:

 (i) $\dN_j (x) >0$ for $j \not\equiv 0~(\bmod~2r).$
 
 (ii) $\dN_{2kr}(x) =0$ for $k=1, 2, 3, ...$

\noindent These conditions are equivalent to requiring  
$\Gamma_{2r} \subset  \ddL$, and that all $x \in \Gamma_{2r}$
have the same balance-set $Z(x) := 
\{ k: D_{k}(x) =0\} = 2r\NN.$

\paragraph{Claim 1.}   {\em All members  $x \in \Gamma_{2r}$ 
have  local level sets  of positive
Hausdorff dimension,  satisfying
$$
 \dim_{H}(L_x^{loc}  )\ge \frac{\log 2}{\log (2^{2r})}= \frac{1}{2r}.
$$}

Claim 1 will follow from the spacing of the balance points
being an arithmetic progression.  This makes  each local level set $L_x^{loc}$
 a Cantor-like set, which has a standard  tree construction covered by $2^k$ intervals
of length $2^{-2rk}$, so that $\dim_{H}(L_x^{\loc}) = \frac{1}{2r}$.
Note that the particular subintervals are chosen differently at each step, so that this
is generally not a self-similar construction, but the Hausdorff dimension lower bound proof 
for the Cantor set  given in Falconer \cite[Sect. 2.3]{Fa03} remains
valid here, establishing 
Claim 1.  \medskip

Claim 1 shows that $\Gamma_{2r} \subset \Lambda_{H}^L$, so that
$$
\bigcup_{r=1}^{\infty} \Gamma_{2r} \subset \Lambda_{H}^L.
$$
To complete the proof it suffices to show the sets $\Gamma_{2r}$ each have positive
Hausdorff dimension, which approaches $1$ as $r \to \infty.$

\paragraph{Claim 2.}{\em  For all large enough $r$, the set $\Gamma_{2r}$ has Hausdorff dimension
greater than $1 - \frac{2 \log r}{r}$.}\medskip

Claim 2  follows by observing that $\Gamma_{2r}$ is itself a 
self-similar Cantor set in
which each block of $2r$ symbols is drawn from the set
$$
X_{2r} := \{ x= \frac{m}{2^{2r}}= 0.b_1 b_2... b_{2r}: ~~ \dN_j(x) > 0 ~\mbox{for} ~ 1\le j< 2r,~~
\dN_{2r}(x) =0\},
$$
whose Hausdorff dimension is computable by the method of Falconer \cite[Sect. 2.3]{Fa03}.
It is well known that 
$$
|X_{2r}| = C_r = \frac{1}{r+1} \left( {{2r}\atop{r}} \right),
$$
is a Catalan number. Thus we obtain
$$
\dim_{H} (\Gamma_{2r}) = \frac{\log C_r}{\log 2^{2r}}= \frac{\log C_r}{2r \log 2}.
$$
However it is well  known that $C_r = \frac{2^{2r}}{ \pi r^{\frac{3}{2}}}\left(1+o(1)\right),$ 
as the integer $r \to \infty$. We
conclude that for large enough $r$ there holds
$$
\dim_{H}(\Gamma_{2r}) > 1 - \frac{2 \log r}{r}.
$$
Claim 2 now  follows, and the theorem is proved.

Note that Claims 1 and 2 together imply that the local level set dimension spectrum
function $f_{\tau}^{\ast}(\alpha)$ defined in Section \ref{sec12} satisfies,
$f_{\tau}^{\ast}(\frac{1}{2r}) > 1- \frac{2 \log r}{r}$ for all sufficiently large integers $r$.
$~~~\Box$.  \medskip

\paragraph{Proof of Theorem \ref{th13a}.}
We have
 $\Lambda_{H}^{L} \subset \Gamma_{H}^{abs}$.   Theorem \ref{th51}  then gives
$\dim_{H}(\Gamma_{H}^{abs}) \ge \dim_{H}(\Lambda_{H}^L) = 1.$
yielding
$\dim_{H}( \Gamma_{H}^{abs}) =1$. $~~~\Box$.

%
%
%
%

\section{Level Sets of  Positive Hausdorff Dimension: Ordinate View}\label{sec7}
\setcounter{equation}{0}

Our object  is  to prove Theorem \ref{th13},
which asserts that the set of ordinate levels $y$ having $\dim_{H}(L(y)) >0$ has
Hausdorff dimension $1$.
We use the result on abscissas proved in the last section (Theorem~\ref{th51}),
together with a property showing that the Takagi function restricted to
certain small domains in $[0,1]$ is quite well behaved, i.e. it is bi-Lipschitz map.
This allows the transfer of Hausdorff dimension lower bounds from the
abscissa case treated in Sect. \ref{sec6}.


\subsection{Bi-Lipschitz property of Takagi function on $\Gamma_{2r}$} \label{sec71}

The following proof of the bi-Lipschitz property incorporates  a substantial
simplification of the authors' original argument, due to P. Allaart.
%
%

\begin{theorem}\label{th52}
Let $Z(x)= \{ j: \, D_j(x)=0 \}$ be the balance-set of $x \in [0, 1]$. 
 For $r \ge 1$ the Takagi function $\tau(x)$ restricted to the (compact) domain
 $$
 \Gamma_{2r} := \{ x \in \ddL:  ~Z(x)= 2r\NN \},
 $$
 is strictly increasing and is a bi-Lipschitz map.
   \end{theorem}

\paragraph{Proof.}
We will show  that if  $x_1 < x_2$ are both in $\Gamma_{2r}$ then
we have the bi-Lipschitz estimates
\beql{711}
2^{2r}  (x_2 - x_1)         \ge  \tau(x_2) - \tau(x_1) \ge   \frac{1}{2^{2r}} ( x_2 - x_1).
\eeq
The lower bound shows that the Takagi function restricted to $\Gamma_{2r}$ is
strictly increasing. 

To prove \eqn{711}, we first reduce to the case that $x_1$ and $x_2$ have
binary expansions that disagree somewhere in the first $2r$ digits.
If they disagree first between the $2kr+1$ and $2(k+1)r$ digits,
then we can write $x_i= x_0 + \frac{w_i}{2^{2kr}}$ with $0 \le w_i \le 1$,
where $x_0=0.b_1 b_2 \cdots b_{2kr}$ is the initial block of $2kr$ digits
where they agree, and $w_1< w_2$.
Note that $w_1, w_2 \in \Lambda_{2r}$ by the self-similar definition of $\Lambda_{2r}$,
and $w_1, w_2$ disagree somewhere in their first $2r$ digits.
Now 
$$
x_2- x_1 = \frac{1}{2^{2kr}} \left(w_2 -w_1\right)
$$
and, since $D_{2kr}( x_1) = D_{2kr}(x_2)=0$, Lemma \ref{le25} gives
$$
\tau(x_2) - \tau(x_1) = \frac{1}{2^{2kr}}\left( \tau(w_2) - \tau(w_1) \right).
$$
Thus it suffices to prove \eqn{711} for $w_1$ and $w_2$.

By definition all  $x =0.b_1 b_2 ...\in \Lambda_{2r}$ 
have $D_m(x) >0$ for all $m \ne 2kr$, while  $D_{2kr}(x) =0$ for all  $k \ge 1$. This forces
 $b_{2kr+1}=b_{2kr+2}=0$
and $b_{2(k+1)r-1}= b_{2(k+1)r- 2}=1$, 
for all $k \ge 0$. We now suppose
$x_1= 0. b_1 b_2 ... $ and $x_2 = 0.b_1^{'} b_2^{'}...$ disagree somewhere
in their first $2r$ digits, say at the $n$-th digit, where $3 \le n \le 2r-2.$
 Now define the dyadic rationals 
$$\tilde{x}_1:= 0. b_1 b_2 ... b_{2r-2} 1\,0^{\infty} = \frac{k_1}{2^{2r-1}},$$
$$\tilde{x}_2 := 0.b_1^{'} b_2^{'} ... b_{2r-2}^{'} 1 \,0^{\infty}= \frac{k_2}{2^{2r-1}},$$
which truncate the $x_i$ at their $(2r-1)$-st digits, and the dyadic rationals 
$$\bar{x}_1:= 0. b_1 b_2 ... b_{2r-2} 11\,0^{\infty} = \frac{2k_1+1}{2^{2r}},$$
$$\bar{x}_2 := 0.b_1^{'} b_2^{'} ... b_{2r-2}^{'} 11\, 0^{\infty}= \frac{2k_2+1}{2^{2r}},$$
which truncate the $x_i$ at their $2r$-th digits.
 We note that $D_{2r-1}(\tilde{x}_i) =1$ 
and $D_{2r}(\bar{x}_i)=0$ for $i=1, 2$.  We can now write
$$
x_i = \bar{x}_i + \frac{ x_i[2r]}{2^{2r}},~~~~\mbox{for} ~~~ i=1, 2,
$$
in which $x_1[2r]= 0.b_{2r+1} b_{2r+2} \cdots $ and similarly for $x_2[2r]$.
Since each $x_{i}[2r] =0.00...$ and each $\bar{x}_i = 0. 00...$ 
we deduce $|x_2[2r] - x_1[2r]| \le \frac{1}{4},$
and
$$
\frac{1}{4} \ge x_2 - x_1 \ge ( \bar{x}_2 - \bar{x_1})- \frac{1}{4} \cdot \frac{1}{2^{2r}}.
$$
Now $\bar{x}_2 - \bar{x}_1 \ge \frac{1}{2^{2r-1}}$ so we obtain
\beql{712}
\frac{1}{4} \ge x_2 - x_1 \ge \frac{7}{8} \cdot \frac{1}{2^{2r-1}}.
\eeq

Next we estimate the difference between $\tau(x_2)$ and $\tau(x_1)$.
For the lower bound, we will relate it to the difference between
$\tau( \bar{x}_2)$ and $\tau(\bar{x}_1)$. Using Lemma \ref{le25}
and $D_{2r}(\bar{x}_i)=0$ for $i=1, 2$, we obtain
$$
\tau(x_i) = \tau( \bar{x}_i) + \frac{ \tau( x_i[2r])}{2^{2r}}, ~~~\mbox{for}~~~i=1,2.
$$
This yields, using $0 \le \tau(x_i[2r]) \le \frac{2}{3}$, 
\beql{713}
\tau( x_2) - \tau(x_1) \ge \left(\tau( \bar{x}_2) - \tau(\bar{x}_1) \right) -\frac{2}{3} \cdot \frac{1}{2^{2r}}.
\eeq
To estimate the right side, we will relate these quantities to $\tau( \tilde{x}_i)$.
Using $\bar{x}_i = \tilde{x}_i + \frac{1}{2^{2r}}$, a calculation using Lemma \ref{le25} gives
$$
\tau( \bar{x}_i) = \tau (\tilde{x}_i) + \frac{1}{2^{2r-1}}, ~~~\mbox{for}~~~i=1,2.
$$
Thus \eqn{713} becomes
\beql{714}
\tau( x_2) - \tau(x_1) \ge 
\left(\tau( \tilde{x}_2) - \tau(\tilde{x}_1) \right) -\frac{2}{3} \cdot \frac{1}{2^{2r}}.
\eeq
Now comes the key observation. Both $\tilde{x}_1, \tilde{x}_2 \in \frac{1}{2} \ddL$,
since $D_m( \tilde{x}_i) \ge 1$ for $1 \le i \le 2r-1$. By hypothesis
$x_1$ and $x_2$ disagree at the $n$-th digit, with $n \le 2r-2$, 
at which  $b_n=0, b_n^{'}=1$, 
and we define the intermediate dyadic rational
$$
x_0 := 0.b_1^{'} b_2^{'}\cdots b_n^{'} \, 0^{\infty} = \frac{k_0}{2^{n}}.
$$
We have $x_0 \in \frac{1}{2} \ddL$ because its binary expansion is a
prefix of that of $x_2$. Furthermore $x_0 < x_2$ because $x_0$ and $x_2$
disagree at the $(2r-1)$-st digit. Thus $x_1 < x_0 < x_2$ and
we may now apply Theorem \ref{th33}(2) to conclude the strict inequality
$\tau(\tilde{x}_1) < \tau(\tilde{x}_2).$
 But $\tau(\tilde{x}_i)$ ($i=1,2$) are both
dyadic rationals with denominator at most $2^{2r-1}$, so we deduce that
$$
\tau( \tilde{x}_2) - \tau( \tilde{x}_1) \ge \frac{1}{2^{2r-1}}.
$$
Substituting this inequality into \eqn{714} yields the lower bound
\beql{715}
\tau( x_2) - \tau(x_1) \ge\frac{4}{3} \cdot \frac{1}{2^{2r}}.
\eeq
This proves the function $\tau(x)$ is strictly increasing on $\Gamma_{2r}$,
and also   gives the desired bi-Lipschitz estimate. 

Combining \eqn{715} with the trivial upper bound estimate 
$\tau(x_2) - \tau(x_1) \le \frac{2}{3},$
and  with \eqn{712} yields
 $$
2^{2r-1} (x_2- x_1) \ge \tau( x_2) - \tau(x_1) \ge \frac{1}{2^{2r-2}} (x_2 -x_1).
 $$
This implies \eqn{711}. $~~~\Box$

\paragraph{Remark.} 
 The Takagi function  $\tau(x)$ is not a  Lipschitz map on its
 full domain $[0,1]$, nor is it  a Lipschitz function even when restricted to the  domain $\ddL$.
 This is because it has arbitrarily steep slopes on $\ddL$, as is implicit in the singular function property. 


\subsection{Hausdorff dimension of $\Gamma_H^{ord}$}\label{sec72}

To conclude the paper we prove Theorem \ref{th13}.

\paragraph{Proof of Theorem~\ref{th13}.} 
Let 
$$
\Gamma_H^{ord} := \{ y: ~0\le y \le \frac{2}{3} ~\mbox{with}~\dim_{H} L(y) > 0\}.
$$
It is well known that bi-Lipschitz maps preserve Hausdorff dimension.
By Theorem~\ref{th52}  the bi-Lipschitz property  holds for the Takagi function $\tau$ 
restricted to the compact domain $\Gamma_{2r}$. The range of this map is
$$
\tilde{\Gamma}_{2r}:= \{ y: y= \tau(x), ~x \in \Gamma_{2r} \},
$$
which therefore satisfies 
$$
\dim_{H}( \tilde{\Gamma}_{2r}) = \dim_{H}(\Gamma_{2r}) \ge  1 - \frac{2 \log r}{2r}. 
$$
for large enough $r$, as shown in Claim 2
of the proof of Theorem~\ref{th51}. 

Next, 
Claim 1 of that proof shows that every level $y \in \Gamma_{2r}$
has
$$
\dim_{H}( L(y) ) \ge \dim_{H} (L_x^{\loc}) \ge \frac{1}{2r}.
$$
Combining these inequalities  shows that the Takagi dimension spectrum function satisfies
$$
f_{\tau}( \frac{1}{2r} ) \ge 1- \frac{2\log r}{2r}
$$
for all sufficiently large $r$.
Thus we have 
$\tilde{\Gamma}_{2r} \subset \Gamma_{H}^{ord}.$
and
$$
\dim_{H} (\Gamma_H^{ord}) \ge 1- \frac{2\log r}{r}.
$$
Letting $r \to \infty$ gives $\dim_{H}(\Gamma_H^{ord}) \ge 1$, establishing the result.
$~~\Box$

%
%
%
%

\noindent 
Jeffrey C. Lagarias \\
Dept. of Mathematics\\
The University of Michigan \\
Ann Arbor, MI 48109-1043\\
\noindent email: {\tt lagarias@umich.edu}\\

\noindent Zachary Maddock \\
Dept. of Mathematics\\
Columbia University \\
New York, NY 10027\\
\noindent email: {\tt maddockz@math.columbia.edu}

\end{document}